\documentclass[12pt,reqno]{amsart}

\usepackage{amsmath}
\usepackage{amssymb}
\usepackage{amsthm}
\usepackage[english]{babel}

\newtheorem{theorem}{Theorem}

\newcommand{\beqa}{\begin{eqnarray}}
\newcommand{\beqan}{\begin{eqnarray*}}
\newcommand{\eeqa}{\end{eqnarray}}
\newcommand{\eeqan}{\end{eqnarray*}}
\def\beq#1\eeq{\begin{equation}#1\end{equation}}

\def\P{\mathbf P }

 \def\na{\,\, {\raise.4pt\hbox{$\shortmid$}}{\hskip-2.0pt\to}\, \, }

\def\={\overset{ \text{\rm def} }=}

\def\ffrac{\frac}
\def\4{\kern1pt}

\def\bgr#1{\4\bigr#1}
\def\bgl#1{\bigl#1\4}
\def\BR{\bf R}

\begin{document}

\title[closeness of convolutions on convex polyhedra]
{Estimates for the closeness of convolutions of probability distributions on convex polyhedra}

\author[F.~G\"otze]{Friedrich G\"otze}
\author[A.Yu. Zaitsev]{Andrei Yu. Zaitsev}

\email{goetze@math.uni-bielefeld.de}
\address{Fakult\"at f\"ur Mathematik,\newline\indent
Universit\"at Bielefeld, Postfach 100131,\newline\indent D-33501 Bielefeld,
Germany\bigskip}
\email{zaitsev@pdmi.ras.ru}
\address{St.~Petersburg Department of Steklov Mathematical Institute
\newline\indent
Fontanka 27, St.~Petersburg 191023, Russia\newline\indent
and St.Petersburg State University, 7/9 Universitetskaya nab., St. Petersburg,
199034 Russia}

\begin{abstract}
The aim of the present work is to show that the results obtained earlier
on the approximation of distributions of sums of independent summands by
the accompanying compound Poisson laws and the estimates of the proximity of sequential convolutions of multidimensional distributions may be transferred to the estimation of the closeness of convolutions of probability distributions on convex polyhedra.

\end{abstract}



\thanks{The authors were supported by the SPbGU-DFG grant 6.65.37.2017.
The second author was supported by grant RFBR 16-01-00367.}

\maketitle

Let us first introduce some notation. Let $\mathfrak F_d$ denote the set
of probability distributions defined on the Borel $\sigma$-field of
subsets of the Euclidean space~${\BR}^d$ and let \,${\mathcal L}(\xi)\in\mathfrak F_d$ \,be the
distribution of a $d$-dimensional random vector~$\xi$. Let
 $\mathfrak F_d^s\subset \mathfrak F_d$ be the set of symmetric distributions.
For $F\in\mathfrak F_d$, we
denote the corresponding characteristic functions by $\widehat F(t)$, $t\in \mathbf R^d$, and distribution functions
by~$F(x)=F\{(-\infty,x_1]\times\cdots\times(-\infty,x_d]\}$, $x=(x_1,\ldots,x_d)\in{\BR}^d$.
The uniform Kolmogorov distance is defined by
$$ \rho (G,H)=\sup_{x\in \mathbf R^d} \;\bgl|G(x)-H(x)\bigr|,\quad G,H\in\mathfrak F_d. $$ By the symbols~$c$ and~$c(\,\cdot\,)$ we denote (generally speaking various) positive absolute constants and quantities depending only on the arguments in brackets.
For $0\le\alpha\le2$,
we denote
$$
\mathfrak F_d^{(\alpha)}=\Big\{F\in \mathfrak F_d^s:\widehat F(t)\ge-1+\alpha,
\ \hbox{for all }t\in \mathbf R^d \Big\},\quad \mathfrak F_d^{+}=\mathfrak F_d^{(1)}.
$$

Products and powers of measures are understood in the
convolution sense: \;${GH=G*H}$, \,$H^m=H^{m*}$, \,$H^0=E=E_0$,
where $E_x$ \,is the distribution concentrated at a point $x\in
{\BR}^d $. A natural approximating infinitely divisible distribution
for $\prod_{i=1}^nF_i$ is the accompanying compound Poisson
distribution $\prod_{i=1}^n e(F_i)$, where $$
e(H)=e^{-1}\sum_{k=0}^\infty
 \ffrac{H^k}{k!},\quad H\in\mathfrak F_d,
$$ and, more generally,
\begin{equation}  e(\alpha\4 H)=e^{-\alpha}\sum_{k=0}^\infty
 \ffrac{\alpha^k\4H^k}{k!}, \quad \alpha>0.  \label{00}\end{equation}
It is well-known that the distribution $e(\alpha\4 H)$ is infinitely divisible.

 Arak \cite{A80} showed that, if $F$ is a symmetric one-dimensional distribution with a nonnegative
characteristic function for all $t \in \mathbf R$, then
\begin{equation}
\rho(F^n,e(n\4F)) \le c\,n^{- 1 },
\label{a97}\end{equation}
He introduced and used the so-called method of triangular functions (see \cite[Chapter~3, Sections~2--4]{2}).

 Zaitsev \cite{z83} applied the methods which were used by Arak while proving inequality \eqref{a97} (see \cite[Chapter~5, Sections~2, 5--7]{2}). Later, he managed to modify these methods, adapting them to the multidimensional case (see \cite{z87}--\cite{z92}). In particular,
in \cite{z892}, a multidimensional analogue of inequality \eqref{a97} was obtained.

Using the method of triangular functions and its generalizations, several bounds of the type \begin{equation}\rho(G,H)\le c(d)\, \varepsilon\label{b97}\end{equation} were obtained, where $0<\varepsilon<1 $ is small, $G, H\in \mathfrak F_d$, and the inequalities
\begin{equation} \sup_{t\in\mathbf{\mathbf R^d}}\;  \bigl|\widehat G(t)-\widehat H(t)\bigr|\le c\,\varepsilon \label{c97}\end{equation} are valid (see a discussion for $d=1$ in \cite[Chapter~3, Section~3]{2}). Note that, in the general case, \eqref{c97} does not imply~\eqref{b97}.

Inequality~\eqref{b97} is equivalent to the validity of the inequality \begin{equation} \bigl|G\{X\}-H\{X\}\bigr|\le c(d)\,\varepsilon \label{g97}\end{equation}
for all sets $X$ of the form
\begin{equation} X=\big\{x\in\mathbf R^d:a_j\le\langle x,e_j\rangle\le b_j, \ j=1,\ldots, d\big\}, \label{g9}\end{equation}
where $e_j\in \mathbf R^d$ are the vectors of the standard Euclidean basis, $-\infty\le a_j\le b_j\le \infty$, $j=1,\ldots, d$.

For $m\in\mathbf N$ we denote by $\mathfrak X_m$ the collection of convex polyhedra $X \subset \mathbf R^d$
representable in the form
$$
X=\big\{x\in\mathbf R^d:a_j\le\langle x,t_j\rangle\le b_j, \ j=1,\ldots, m\big\},
$$
where $t_j\in \mathbf R^d$, $-\infty\le a_j\le b_j\le \infty$, $j=1,\ldots, m$,
and, for $H={\mathcal L}(\xi)\in\mathfrak F_d$, $X\in\mathfrak X_m$,
$$
q(H, X)=\inf_{t\in \mathbf R^d,\\ \|t\|=1} Q\big(\mathcal L\big(\langle \xi,t\rangle\big),
\lambda\big\{\{\langle x,t\rangle:x\in X\}\big\}\big),
$$
where $\lambda\{\,\cdot\,\}$ is the Lebesgue measure
 and \;${Q(F,\,b)=\sup_x
F\bgl\{[x,\,x+b]\bgr\}}$ is the
concentration function of $F\in\mathfrak F_1$.
Define a distance over all polyhedra $X\in\mathfrak X_m$ as
$$ \rho_m (G,H)=\sup_{X\in\mathfrak X_m} \bigl|G\{X\}-H\{X\}\bigr|. $$

 In~\cite{z92}, Zaitsev has proved several non-uniform inequalities
of the form $| G \{X\} - H \{X\}|\le c(m)\,\varepsilon\,\beta(G,H,X)+o(\varepsilon)$ containing a factor $\beta(G,H,X)$,
depending on the compared distributions
and on the set $X\in\mathfrak X_m$, which satisfies the inequality $\beta(G,H,X)\le c(m)$
and may turn out to be small if the
polyhedron $X$ is sufficiently small in a certain sense.
Theorems~\ref{Th0} and~\ref{Th1} have been obtained in~\cite{z92} as consequences of the corresponding results of~\cite{z88}--\cite{z892} which were proved for the case when
$X$ is a parallelepiped \eqref{g9} with faces parallel to the coordinate axes. The aim of the present paper is to formulate and to discuss similar bounds for the quantities $\rho_m (G,H)$ and $\bigl|G\{X\}-H\{X\}\bigr|$,
$X\in\mathfrak X_m$, where $G, H\in \mathfrak F_d$ are certain convolutions of probability distributions which have not been studied in the literature before.

\begin{theorem}\label{Th0} Let  $F\in\mathfrak F_d^{(\alpha)}$,
$0\le\alpha\le2$,
$m,n\in\mathbf N$, $X\in\mathfrak X_m$,  $D=e(nF)$,  $q_1= q(D, X)$. Then
\begin{equation} \bigl| (F^n) \{X\} - D \{X\}\big|  \le c(m)\,\Big(n^{-1} q_1^{1/5}
\big(\bigl|\log q_1\bigr|+1\big)^{(17m+24)/5}+
\exp(-n\alpha+ c\,m\,\log^3 n)\Big).  \label{999}\end{equation}
and, moreover,
\begin{equation} \bigl| (F^n) \{X\} - (F^{n+1}) \{X\}\big|  \le c(m)\,\Big(n^{-1} q_1^{1/3}
\big(\bigl|\log q_1\bigr|+1\big)^{3m+2}+
\exp(-n\alpha+ c\,m\,\log^3 n)\Big).
  \label{997}\end{equation}
Hence,
\begin{equation} \max\big\{ \rho_m(F^n,e(nF)),  \rho_m(F^n,F^{n+1})\big\} \le c(m)\,\Big(n^{-1} +
\exp(-n\alpha+ c\,m\,\log^3 n)\Big).  \label{99}\end{equation}
\end{theorem}

\begin{theorem}\label{Th1} Assume that the
 distributions $G_i\in\mathfrak F_d$ are represented as
 \begin{equation}
G_i=(1-p_i)\,E+p_i\,V_i, \label{11h}\end{equation}
where $V_i\in\mathfrak F_d$ are arbitrary distributions, $0\le p_i\le p=\max_jp_j$,
$$m\in\mathbf N,\quad
X\in\mathfrak X_m,\quad G=\prod_{i=1}^nG_i,\quad D=\prod_{i=1}^ne(G_i),   \quad q_2= q(D_0, X),$$
where $D_0$ are the $d$-variate infinitely divisible distribution with characteristic
function
$$
\widehat D_0(t)=\prod_{i=1}^n\exp\big(-p_i(1-p_i)\big(1-\hbox{\rm Re}\,\widehat V_i(t)\big)\big), \quad t\in \mathbf R^d.
$$
Then
\begin{equation} \bigl| G \{X\} - D \{X\}\big|  \le c(m)\, q_2^{1/3}
\big(\bigl|\log q_2\bigr|+1\big)^{3m+2}\,p  \label{699}\end{equation}
and, hence,
\begin{equation} \rho_m(G,D)  \le c(m)\, p.  \label{r699}\end{equation}
\end{theorem}

It is easy to see that, under the conditions of Theorems \ref{Th0} and~\ref{Th1}, we have
$0 \le q_j \le 1$, $j = 1, 2$, and,
moreover, the quantities $q_j$ may be small. For example, for a fixed bounded set
$X\in\mathfrak X_m$, the quantity $q_1$ decreases
for $n\to\infty$ not slower than $O(n^{-1/2})$. Thus, inequalities \eqref{999}, \eqref{997} and~\eqref{699} significantly strengthen inequalities \eqref{99} and~\eqref{r699}. At the
same time, there is no reason to expect that inequalities \eqref{999}, \eqref{997} and~\eqref{699} are optimal with respect to the dependence of the right-hand
sides on the parameters $q_1$ and~$q_2$. In particular, already from Arak's results
(see \cite[Theorem 7.1, Chap. V]{2}) it follows that, for $m=1$, $\alpha=1$, inequality~\eqref{999} may be replaced by
\begin{equation} \bigl| (F^n) \{X\} - D \{X\}\big|  \le c\,n^{-1} q_1^{1/3}
\big(\bigl|\log q_1\bigr|+1\big)^{13/3}.  \label{909}\end{equation}

In inequalities \eqref{999} and \eqref{997} we have $\varepsilon=n^{-1}$ and
$$\beta(F^n,D,X)=q_1^{1/5}\big(\bigl|\log q_1\bigr|+1\big)^{(17m+24)/5}, $$
$$\beta(F^n,F^{n+1},X)=q_1^{1/3}\big(\bigl|\log q_1\bigr|+1\big)^{3m+2} $$ respectively.

The proofs of Theorems \ref{Th0} and~\ref{Th1} are based on applications of $m$-variate versions of bounds for the closeness of distributions on the sets of the form~\eqref{g9} with $d=m$. It is important that the $m$-variate vectors with coordinates $\langle \xi,t_j\rangle$, $\langle \eta,t_j\rangle$, $t_j\in\mathbf R^d$, $j=1,\ldots, m$, satisfy the same $m$-dimensional conditions as the random vectors $\xi,\eta\in\mathbf R^d$ with compared $d$-dimensional distributions. For example, if
$F=\mathcal{L}(\xi)\in\mathfrak F_d^{(\alpha)}$,
for some $\alpha$ satisfying $0\le\alpha\le2$, then $\mathcal{L}\big(\langle \xi,t_1\rangle,\ldots,\langle \xi,t_m\rangle\big)\in\mathfrak F_m^{(\alpha)}$.
Similarly, if
$F=\mathcal{L}(\xi)\in\mathfrak F_d^{s}$,
 then $\mathcal{L}\big(\langle \xi,t_1\rangle,\ldots,\langle \xi,t_m\rangle\big)\in\mathfrak F_m^{s}$. Analogous statements hold for $n$ and $(n+1)$-fold convolutions of such distributions and about other distributions involved in the assertions of Theorems \ref{Th0} and~\ref{Th1}. Thus, roughly speaking, from the known estimates of the distance $\rho$ in space $\mathbf R^m$ we derive estimates of the distance $\rho_m$ in space $\mathbf R^d$.

The situation considered in Theorem~\ref{Th1} can be interpreted as a comparison of the sample containing
independent observations of rare events with a Poisson point process
which is obtained after a Poissonization of the initial sample (see~\cite{GZ17},~\cite{14}). 

Indeed, let $X_1, X_2,\dots, X_n$ be independent not
identically distributed elements of a measurable space $(\mathfrak
X,\mathcal S)$ and $f:\mathfrak X\to\mathbf R^m$ be a  Borel mapping.
Assume that the set $\mathfrak X$ is represented as the union of two disjoint measurable sets: $\mathfrak X= \mathfrak
X_1\cup\mathfrak X_2$, with \,$\mathfrak X_1,\,\mathfrak X_2\in \mathcal S$, $\mathfrak
X_1\cap\mathfrak X_2=\varnothing$. We say that the $i$-th  rare event occurs if $X_i\in\mathfrak X_2$. 
Respectively, it does not occur if $X_i\in\mathfrak X_1$.

Assume that $f(x)=0$, 
for all $x\in\mathfrak X_1$, and
$$0\le p_i=\mathbf{P}\big\{X_i\in\mathfrak X_2\big\}
=1-\mathbf{P}\big\{X_i\in\mathfrak X_1\big\}\le p=\max_{1\le i\le
n}p_i. \nonumber$$ Then $\mathcal{L}\big(f(X_i)\big)=(1-p_i)\,E+p_i\,V_i$, where $E, V_i\in\mathfrak F_m$, and $V_i$ 
is conditional distributions of $f(X_i)$ given
 $X_i\in\mathfrak X_2$. In~\cite{GZ17}, it was shown that 
\begin{equation}\rho\bigg(\mathcal{L}\Big(\sum_i f(X_i)\Big),
\mathcal{L}\Big(\sum_k f(Y_k)\Big)\bigg)\le c(m)\,p,\label{99h}\end{equation}
where $Y_k$ are the points of the corresponding Poisson point process. In particular, in the case where $\mathcal{X}=\mathbf R^d$, $\mathcal{X}_1=\{0\}$,
and $f(x)=\big(\langle x,t_1\rangle,\ldots,\langle x,t_m\rangle\big)$, for $x\in\mathbf R^d$,
inequality~\eqref{99h} turns into inequality~\eqref{r699}.
\medskip

In the rest of the paper we will study how
small is the difference between $F^{n +k}$ and $F ^n$. A particular case of this problem is considered in inequality \eqref{999} of Theorem~\ref{Th0}.

In the papers of Zaitsev \cite{z80}, \cite{z81}, \cite{z88} it was shown that
one can obtain sharp 
bounds for the closeness of $F^{n +k}$ and $F ^n$ without any moment conditions. Moreover, if the distribution
$F$ is centered so that all its marginal distributions have zero medians, then
\begin{equation}
\rho(F^n,F^{n+1}) \le c\,d\,n^{- 1/2 },
\label{609}\end{equation}
where $c$ is an absolute constant (which was estimated for $d=1$ in \cite{Mai} and \cite{Gol}: $c\le5.85$).
The proof of this inequality is
relatively simple and is based on classical bounds for concentration functions of
convolutions.
Much more complicated methods are needed to investigate the case of
symmetric distributions $F\in\mathfrak F_d^s$. In
this case inequality~\eqref{609} is valid and
it is optimal
with respect to order in~$n$. But it may be essentially improved in the
case when the characteristic function $\widehat F(t)$ is uniformly separated from $-1$. In
particular,
\begin{equation}
\rho(F^n,F^{n+1}) \le c(d)\,n^{- 1 }.
\label{639}\end{equation}
if $\widehat F(t)\ge0$ for all $t\in \mathbf R^d $. Notice that inequality~\eqref{997} is much more general compared to~\eqref{639}.
Using
this fact for the distribution $F^2$ with symmetric $F$, we obtain the paradoxical
statement that for all natural numbers $n$ and for any symmetric distribution $F$ the
inequalities
\begin{equation}
\rho(F^n,F^{n+1}) \le c\,d\,n^{- 1/2 }\quad\hbox{and}\quad\rho(F^n,F^{n+2}) \le c(d)\,n^{- 1 }
\label{b950}\end{equation}
are valid and they are both optimal with respect to the order in~$n$. Inequalities~\eqref{b950} imply the following Theorem~\ref{Th9}.

\begin{theorem}\label{Th9} Let  $F\in\mathfrak F_d^{s}$,
$k,n\in\mathbf N$. Then
\begin{equation}  \rho(F^n,F^{n+2k}) \le c(d)\,k\,n^{-1} .  \label{b90}\end{equation}
\begin{equation}  \rho(F^n,F^{n+2k+1}) \le  c\,d\,n^{- 1/2 }+c(d)\,k\,n^{-1} .  \label{c90}\end{equation}In particular,
\begin{equation} \sup_{k\le\sqrt{n}} \rho(F^n,F^{n+k}) \le  c(d)\,n^{- 1/2 } .  \label{c907}\end{equation}
\end{theorem}

It is evident that knowledge about the closeness of $F^{n +k}$ and $F ^n$ is useful for
studying distributions of the form
$$
G=\sum_{s=0}^\infty p_s\,F^s, \quad 0\le p_s\le1,\quad\sum_{s=0}^\infty p_s=1.
$$
A result in this direction is given in our Theorem~\ref{Th888}.
In particular, using \eqref{00} and the bounds for the closeness of $F^{n +k}$ and $F ^n$, Zaitsev~\cite{z88} proved the following Theorem~\ref{Th33}.

\begin{theorem}\label{Th33} Let  $F\in\mathfrak F_d^{s}$,
$n\in\mathbf N$. Then
\begin{equation}  \rho(F^n,e(nF)) \le c(d)\,n^{-1/2} .  \label{990a}\end{equation}
\end{theorem}

A one-dimensional version of Theorem~\ref{Th33} was proved somewhat earlier in~\cite{z81}.\medskip

 It is evident that if a distribution $F\in\mathfrak F_d$
is concentrated on a hyperplane which does not contain zero and is
orthogonal to one of coordinate axes then
$\rho(F^n,F^{n+k}) =1$ for any $n, k\in\mathbf N$.
In particular, this is true in the case where $F=E_a$, $a\in \mathbf R^d$, $a\ne0$.
On the other
hand, if all distributions $F^{(j)}\in \mathfrak F_1$, $j=1,\ldots, d$, of coordinates of the vector $\xi$ with $\mathcal{L}(\xi)=F$ are either non-degenerate
or equal to~$E\in \mathfrak F_1$,  then, as is shown in Zaitsev~\cite{z80}, $\rho(F^n,F^{n+1})\to0$
as $n\to\infty$ and, moreover,
\begin{equation}
\rho(F^n,F^{n+1}) \le \frac{c(F)}{\sqrt n},\quad\hbox{for all }n\in\mathbf N.
\label{511h}\end{equation}

Let $F\in \mathfrak F_1$ be a one-dimensional lattice symmetric distribution concentrated on the set of odd numbers. Then the distributions $F^n$, $n = 1, 2,\ldots$, are concentrated either on the set of odd numbers or on the set of even ones according to the parity of the number~$n$. Therefore, $\rho(F^n,F^{n+1}) \ge Q(F^n,0)/2$.
For many distributions, e.g. for $F=E_{-1}/2+E_1/2$, the concentration function $Q(F^n,0)$ behaves as
$c(F)\,n^{-1/2}$ as $n\to\infty$. This indicates that the rate of decrease
with respect to $n$ of the right-hand side of
\eqref{511h} cannot be increased without additional assumptions.

It is easy to show that the distribution $F\in\mathfrak{F}_1^s $ is concentrated on the set of odd numbers if and only if its characteristic function $\widehat F(t)$ is equal to $-1$ at the points $t=(2k+1)\pi$ for all $k\in\mathbf Z$. For example, let $\widehat F(t)=\cos t$, for $F=E_{-1}/2+E_1/2$. Inequality~\eqref{99} of Theorem~\ref{Th0} says that the separation from $-1$ of the characteristic function of a distribution $F\in\mathfrak F_d^s$ leads
to more quick decay of $\rho(F^n,F^{n+1})$ than
the inequality \eqref{511h} is able to provide.

 Similarly to the proof of Theorems~\ref{Th0} and~\ref{Th1}, we can use inequalities \eqref{b90}--\eqref{511h} in order to obtain the corresponding analogues of \eqref{b90}--\eqref{511h} for the closeness of convolutions of $d$-dimensional distributions on the convex polyhedra $X\in\mathfrak X_m$. The following Theorems~\ref{Th89}--\ref{Th888} are the main results of the present paper.

 \begin{theorem}\label{Th89}
Let  $F\in\mathfrak F_d^{s}$,
$k,m,n\in\mathbf N$. Then
\begin{equation}  \rho_m(F^n,e(nF)) \le c(m)\,n^{-1/2} ,  \label{990}\end{equation}
\begin{equation}    \rho_m(F^n,F^{n+2k})\le c(m)\,k\,n^{-1} ,  \label{b95}\end{equation}
\begin{equation}  \rho_m(F^n,F^{n+2k+1})\le  c\,m\,n^{- 1/2 }+c(m)\,k\,n^{-1} .  \label{c903}\end{equation}In particular,
\begin{equation}  \sup_{k\le\sqrt{n}} \rho_m(F^n,F^{n+k})\le  c(m)\,n^{- 1/2 } .  \label{c904}\end{equation}
\end{theorem}

 For $m\in\mathbf N$, $t_1,\ldots,t_m\in\mathbf R^d$, we denote by $\mathfrak X(t_1,\ldots,t_m)$ the collection of convex polyhedra $X \subset \mathbf R^d$
representable in the form
$$
X=\big\{x\in\mathbf R^d:a_j\le\langle x,t_j\rangle\le b_j, \ j=1,\ldots, m\big\}.
$$
 Clearly,
 $$
\mathfrak X_m=\bigcup_{t_1,\ldots,t_m}\mathfrak X(t_1,\ldots,t_m).
$$
The following Theorem~\ref{Th88} is  a consequence of inequality~\eqref{511h}.

\begin{theorem}\label{Th88} Let  $F=\mathcal{L}(\xi)\in\mathfrak F_d$, $m\in\mathbf N$, $t_1,\ldots,t_m\in\mathbf R^d$, and  all distributions of the random variables $\langle \xi,t_j\rangle$, $j=1,\ldots, m$, are either non-degenerate
or equal to~$E\in \mathfrak F_1$. Then,
for all
$X\in\mathfrak X(t_1,\ldots,t_m)$,
\begin{equation} \bigl| (F^n) \{X\} - (F^{n+1}) \{X\}\big|  \le c(F, t_1,\ldots,t_m)\,n^{-1/2} .  \label{333}\end{equation}
\end{theorem}

Thus, we have the alternative: the left-hand side of \eqref{333} is equal to one or decreases at least as $O(n^{-1/2})$.

The quantity $c(F, t_1,\ldots,t_m)$ can be larger than any absolute constant. For example, if $F=F_n\in\mathfrak F_1$ depends on~$n$ and $F_n\big\{[n,n+1)\big\}=1$, then $\rho_1(F_n^n, F_n^{n+1})=1$.

Note, however, that there is a difference between Theorems \ref{Th0}--\ref{Th1} and Theorems \ref{Th89}--\ref{Th88}. In Theorems \ref{Th0}--\ref{Th1}, the bounds are non-uniform. They include factors $\beta(\,\cdot\,,\,\cdot\,,X)$ which depend on $q_j$, $j=1,2$, and may be small for small sets $X\in\mathfrak X_m$. The bounds of Theorems \ref{Th89}--\ref{Th88} cannot be improved even if we compare the probabilities to hit the set containing the unique point 0 only, where $d=1$ and $F=E_{-1}/2+E_1/2$. An exception is inequality~\eqref{b95}. Applying inequality~\eqref{997} to the distribution $F^2\in\mathfrak F_d^{+}$, it is easy to show that
\begin{equation} \bigl| (F^n) \{X\} - (F^{n+2k}) \{X\}\big|  \le c(m)\,k\,\Big(n^{-1} q_3^{1/3}
\big(\bigl|\log q_3\bigr|+1\big)^{3m+2}+
\exp(-n+ c\,m\,\log^3 n)\Big),
  \label{997a}\end{equation}
for any $X\in\mathfrak X_m$, where  $q_3= q(e(n_0F^2), X)$ and $n_0$ is the maximal integer which is less or equal to $n/2$.

In conclusion, we formulate a result on the closeness of the distributions of sums of random number of independent identically distributed random vectors, which follows from Theorem~\ref{Th89}.
For the distance $\rho(\,\cdot\,,\,\cdot\,)$ this result is contained in \cite[Theorem 1.3]{z88}.

Let $\xi_1, \xi_2, \ldots$ be independent identically distributed random vectors with common distribution $F\in\mathfrak F_d$ and let $(\mu, \nu)\in\mathbf{Z}^2$ be a two-dimensional random vector with integer non-negative coordinates, independent of the sequence $\{\xi_j\}_{j=1}^\infty$. Denote
\begin{equation}
U=\mathcal{L}(\mu), \enskip V=\mathcal{L}(\nu), \enskip G=\mathcal{L}(\xi_1+\cdots+\xi_\mu)
, \enskip H=\mathcal{L}(\xi_1+\cdots+\xi_\nu).
  \label{007a}\end{equation}
Then it is well known that then
\begin{equation}
G=\sum_{k=0}^\infty\P\{\mu=k\}\,F^k,\quad H=\sum_{k=0}^\infty\P\{\nu=k\}\,F^k.
  \label{008a}\end{equation}

\begin{theorem}\label{Th888}  If $F\in\mathfrak F_d^s$, then
\begin{equation}
\rho_m(G,H) \le\inf{\mathbf E}\,\min\Big\{\frac{c\4m}{\sqrt{\nu+1}}+c(m)\,\frac{\left|\mu-\nu\right|}{\nu+1}, 1\Big\},
  \label{708a}\end{equation}
and if $F\in\mathfrak F_d^+$, then
\begin{equation}
\rho_m(G,H) \le\inf{\mathbf E}\,\min\Big\{c(m)\,\frac{\left|\mu-\nu\right|}{\nu+1}, 1\Big\}.
  \label{608a}\end{equation}
Here, the infimum is taken over all possible two-dimensional distributions $\mathcal{L}((\mu, \nu))\in\mathfrak F_2$ such that  $\mathcal{L}(\mu)=U$, $\mathcal{L}(\nu)=V$.
\end{theorem}

An interesting problem is to expand our inequalities to arbitrary convex sets~$X$. This does not follow from Theorems \ref{Th0}, \ref{Th1}, and \ref{Th89}--\ref{Th888}, since the constants $c(m)$ depend on~$m$.


\begin{thebibliography}{20}

\bibitem{A80} T. V. Arak, {\it  On the approximation of $n$-fold convolutions of distributions,
 having a non-negative characteristic function with accompanying laws}.
--- Theory Probab. Appl.   \textbf{25}, No. 2 (1980), 221--243.

\bibitem{2}
T. V. Arak, A. Yu. Zaitsev, {\it  Uniform limit theorems for sums of independent random variables}. ---
Proc. Steklov Inst. Math., \textbf{174} (1988), 222 p.


\bibitem{GZ17}
F. G\"otze, A. Yu. Zaitsev, \textit{Rare events and Poisson point
processes}.
---  Zap. Nauchn. Semin. POMI {\bf466}  (2017),  109--119 (in Russian). English
transl. in J. Math. Sci. (N. Y.), (2019), to appear.

\bibitem{Gol} Ia. S. 
Golikova, \textit{On improvement of the estimate of the distance between
sequential sums of independent random variables}.
---  Zap. Nauchn. Semin. POMI {\bf474}  (2018),  118--123 (in Russian). English
transl. in J. Math. Sci. (N. Y.), (2019), to appear.

\bibitem{Mai}
E. L. Maistrenko,   {\it  Estimate for the absolute constant in an inequality for the uniform distance between distributions of sequential sums of independent random variables}. ---
Zap. Nauchn. Semin. POMI {\bf454} (2016), 216--219 (in Russian). English
transl. in
J. Math. Sci. (N. Y.), {\bf229}, 6  (2018), 741--743.


\bibitem{z80}
A. Yu. Zaitsev,   {\it  The estimation of proximity of distribution of sequential sums of independent identically distributed random vectors}. ---
Zap. Nauchn. Semin. POMI {\bf97} (1980), 83--87 (in Russian). English
transl. in
J. Soviet Math., {\bf24}, 5 (1984), 536--539.

\bibitem{z81}   A. Yu. Zaitsev,   {\it  Some properties of $n$-fold convolutions of distributions}. ---
{Theory Probab. Appl.} {\bf26}, 1 (1981), 148--152.

\bibitem{z83}   A. Yu. Zaitsev,  {\it   On the accuracy of approximation of distributions of sums
of independent random variables---which are nonzero with
a small probability---by means of accompanying laws}. ---
{Theory Probab. Appl.} {\bf28}, 3 (1984), 657--669.


\bibitem{z87}  A. Yu. Zaitsev,   {\it  Multidimensional generalized method of triangular functions}. ---
Zap. Nauchn. Semin. LOMI {\bf 158} (1987), 81--104 (in Russian). English
transl. in J. Soviet Math., {\bf43}, 6 (1988), 2797--2810.

\bibitem{z88}
A. Yu. Zaitsev, {\it  Estimates for the closeness of successive convolutions
of multidimensional symmetric distributions}. ---
Probab. Theory Relat. Fields {\bf79},~2 (1988),   175--200.

\bibitem{z89}
A. Yu. Zaitsev, {\it   Multivariate version of the second Kolmogorov's uniform limit theorem}. ---
Theory Probab. Appl. {\bf34},~1 (1989),   108--128.

\bibitem{z892}  A. Yu. Zaitsev,   {\it  On the approximation of convolutions of multi-dimensional symmetric distributions by accompaning laws}. ---
Zap. Nauchn. Semin. LOMI {\bf 177} (1989), 55--72 (in Russian). English
transl. in J. Soviet Math., {\bf61}, 1 (1992), 1859--1872.

\bibitem{z92}  A. Yu. Zaitsev,   {\it  Certain class of nonuniform estimates in multidimensional limit theorems}. ---
Zap. Nauchn. Semin. POMI {\bf 184} (1990), 92--105 (in Russian). English
transl. in J. Math. Sci. (N. Y.), {\bf68}, 4 (1994), 459--468.


\bibitem{14}  A. Yu. Zaitsev,   {\it  On approximation of the sample by a Poisson point process}. ---
Zap. Nauchn. Semin. POMI {\bf 298} (2003), 111--125 (in Russian). English
transl. in J. Math. Sci. (N. Y.), {\bf128}, 1 (2006), 2556--2563.








\end{thebibliography}
\end{document}